\documentclass[11pt,a4paper]{article}
\usepackage{amsmath, amsfonts, amsthm, latexsym}
\newtheorem{theorem}{Theorem}[section]
\newtheorem{proposition}[theorem]{Proposition}
\newtheorem{corollary}[theorem]{Corollary}

\newtheorem{remark}[theorem]{Remark}

\newtheorem{definition}[theorem]{Definition}
\usepackage{amscd}

\begin{document}
\title{\sc When every multilinear mapping is multiple summing}

\author{Geraldo Botelho\thanks{Supported by CNPq Project
202162/2006-0.}\,\, and Daniel Pellegrino\thanks{Supported by CNPq
Projects 471054/2006-2 and 308084/2006-3.\hfill
}}

\date{}
\maketitle \vspace*{-1.0em}

\begin{abstract}
In this paper we give a systematized treatment to some coincidence
situations for multiple summing multilinear mappings which extend,
generalize and simplify the methods and results obtained thus far.
The application of our general results to the pertinent particular
cases gives several new coincidences as well as easier proofs of
some known results.
\end{abstract}

\section*{Introduction}
\indent \indent Multiple summing multilinear mappings between
Banach spaces have been proved to be a very important and very
useful nonlinear generalization of the ideal of absolutely summing
linear operators (see \cite{acosta, bpgv, defant, collect,
mariodaniel, marceladaniel, davidtese, davidstudia, jmaa, arkiv,
marcelatese}). This class was introduced, independently, by Matos
\cite{collect} (under the terminology {\it fully summing
multilinear mappings}) and Bombal, Per\'{e}z-Garc\'{i}a and
Villanueva \cite{bpgv}. The original methods and deep results due
to Per\'{e}z-Garc\'{i}a \cite{davidtese}, which were a source of
inspiration to us in this paper, have played a crucial role in the
development of the
theory.\\
\indent A coincidence situation for multiple summing mappings is a
situation in which every $n$-linear mapping from $E_1 \times
\cdots \times E_n$ to $F$, where $E_1, \ldots, E_n$ and $F$ are
fixed Banach spaces, is multiple $(q;p_1, \ldots, p_n)$-summing
for some numbers $q, p_1, \ldots, p_n$. It happens that the
condition enjoyed by multiple summing mappings by definition (see
Definition \ref{definition}) is a very restrictive one, so
coincidence situations are supposed to be very rare. Nevertheless,
some situations like that are known (along the paper we will came
through some of them) and in this paper we will prove some more.
Such multilinear concidence theorems are usually proved with the
help of linear coincidence situations. In this paper we give a
unified treatment to this approach, in the sense that we identify
general linear conditions from which multilinear coincidences will
follow. In this fashion we obtain new multilinear coincidence
situations as well as generalizations and simplifications of some
known ones.

\section{Background and notation}
Throughout this paper $n$ is a positive integer, $E_1, \ldots,
E_n, E$ and $F$ will stand for Banach spaces over $\mathbb{K} =
\mathbb{R}$ or $\mathbb{C}$, and $E'$ is the dual of $E$. By
${\cal L}(E_1, \ldots, E_n;F)$ we denote the Banach space of all
continuous $n$-linear mappings from $E_1 \times \cdots \times E_n$
to $F$ with the usual sup norm. If $E_1 = \cdots = E_n = E$, we
write ${\cal L}(^nE;F)$ and if $F = \mathbb{K}$ we simply write
${\cal L}(E_1, \ldots, E_n)$ and ${\cal L}(^nE)$. For the general
theory of multilinear mappings we refer to Dineen \cite{din}.\\
\indent Let $p \geq 1$. By $\ell_p(E)$ we mean the Banach space of
all absolutely $p$-summable sequences $(x_j)_{j=1}^\infty$, $x_j
\in E$ for all $j$, with the norm $\|(x_j)_{j=1}^\infty\|_p =
\left (\sum_{j=1}^\infty \|x_j\|^p \right)^{1/p}$. $\ell_p^w(E)$
denotes the Banach space of all sequences $(x_j)_{j=1}^\infty$,
$x_j \in E$ for all $j$, such that $(\varphi(x_j))_{j=1}^\infty
\in \ell_p$ for every $\varphi \in E'$ with the norm
\[\|(x_j)_{j=1}^\infty\|_{w,p} = \sup\{\|(\varphi(x_j))
_{j=1}^\infty\|_p: \varphi \in E', \|\varphi\| \leq 1\}.\]

\begin{definition}\label{definition}\rm Let $1 \leq p_j \leq q$, $j =1,
 \ldots, n$. An $n$-linear mapping $A \in {\cal L}(E_1,
\ldots, E_n;F)$ is {\it multiple $(q;p_1, \ldots, p_n)$-summing}
if there is a constant $C \geq 0$ such that
$$\left ( \sum_{j_1, \ldots, j_n =1}^{m_1, \ldots, m_n} \left\|
A(x_{j_1}^{(1)}, \ldots, x_{j_n}^{(n)}) \right \|^q \right
)^\frac{1}{q} \leq C \prod_{k=1}^n \left \|
(x_j^{(k)})_{j=1}^{m_k}\right \|_{w,p_j}$$ for every $m_1, \ldots,
m_n \in \mathbb{N}$ and any $x_{j_k}^{(k)} \in E_k$, $j_k = 1,
\ldots, m_k$, $k = 1, \ldots, n$. It is clear that we may assume
$m_1 = \cdots = m_n$. The infimum of the constants $C$ working in
the inequality is denoted by $\pi_{q;p_1, \ldots, p_n}(A)$.
\end{definition}
The subspace $\Pi_{q;p_1, \ldots, p_n}^n(E_1, \ldots, E_n;F)$ of
${\cal L}(E_1, \ldots, E_n;F)$ of all multiple $(q;p_1, \ldots,
p_n)$-summing becomes a Banach space with the norm $\pi_{q;p_1,
\ldots, p_n}(\cdot)$. If $p_1 = \cdots = p_n = p$ we say that $A$
is multiple $(q;p)$-summing and write \linebreak $A \in
\Pi_{q;p}^n(E_1, \ldots, E_n;F)$. The symbols $\Pi_{q;p_1, \ldots,
p_n}^n(^nE;F)$, $\Pi_{q;p}^n(^nE;F)$, \linebreak$\Pi_{q;p_1,
\ldots, p_n}^n(E_1, \ldots, E_n)$, $\Pi_{q;p}^n(E_1, \ldots,
E_n)$, $\Pi_{q;p_1, \ldots, p_n}^n(^nE)$ and $\Pi_{q;p}^n(^nE)$
are defined in the obvious way.

Making $n=1$ we recover the classical ideal of absolutely
$(q;p)$-summing linear operators, for which the reader is referred
to Diestel, Jarchow and Tonge \cite{djt}. For the space of
absolutely $(q;p)$-summing linear operators from $E$ to $F$ we
shall write $\Pi_{q;p}(E;F)$ rather than $\Pi_{q;p}^1(E;F)$.

\begin{remark}\label{remark}\rm Throughout the paper we will obtain
multilinear concidences from linear ones. On the other hand, it
must be clear that multilinear coincidences always imply linear
ones. More precisely, it is not difficult to prove (see the proof
of \cite[Theorem 4.3]{marceladaniel}) that if ${\cal L}(E_1,
\ldots, E_n;F) = \Pi_{q;p_1, \ldots,p_n}^n(E_1, \ldots, E_n;F)$,
then ${\cal L}(E_j;F) = \Pi_{q;p_j}(E_j;F)$, $j = 1, \ldots, n$.
\end{remark}

\section{General results}
Our first result establishes the conditions from which several
(known and new) coincidence theorems will follow.

\begin{theorem}\label{teorema} Let $p, r \in [1,q]$ and let $F$ be a Banach
space. By $B(p,q,r,F)$ we mean the collection of all Banach spaces
$E$ such that
$${\cal L}(E;F) = \Pi_{q;p}(E;F){\it ~and~}{\cal L}(E;\ell_q(F)) = \Pi_{q;r}(E;\ell_q(F)).$$
Then, for every $n \geq 2$,
$${\cal L}(E_1, \ldots, E_n;F) = \Pi_{q;r, \ldots,r,p}^n(E_1, \ldots,
E_n;F)$$ whenever $E_1, \ldots, E_n \in B(p,q,r,F)$.
\end{theorem}

\begin{proof} We proceed by induction on $n$. Case $n=2$: let $E_1, E_2
\in B(p,q,r,F)$. By the open mapping theorem there are constants
$C_1$ and $C_2$ such that
$$\pi_{q;p}(u) \leq C_1\|u\| {\rm ~for~every~}u \in {\cal L}(E_2;F){\it ~and~}  $$
$$\pi_{q;r}(v) \leq C_2\|v\| {\rm ~for~every~}v \in {\cal L}(E_1;\ell_q(F)).$$
Let $A \in {\cal L}(E_1,E_2;F)$. Given two sequences
$(x_j^{(1)})_{j=1}^\infty \in \ell_r^w(E_1)$ and
$(x_j^{(2)})_{j=1}^\infty \in \ell_p^w(E_2)$, fix $m \in
\mathbb{N}$ and consider the continuous linear operator
$$A_1^{(m)} \colon E_1 \longrightarrow \ell_q(F) ~:~ A_1^{(m)}(x) =
(A(x,x_1^{(2)}), \ldots, A(x,x_m^{(2)}),0,0, \ldots). $$ So,
$A_1^{(m)}$ is $(q;r)$-summing and $\pi_{q;r}(A_1^{(m)}) \leq
C_2\|A_1^{(m)}\|$. For each $x \in B_{E_1}$, consider the
continuous linear operator
$$A_x \colon E_2 \longrightarrow F ~:~ A_x(y) = A(x,y). $$
So, $A_x$ is $(q;p)$-summing and $\pi_{q;p}(A_x) \leq C_1\|A_x\|
\leq C_1\|A\|\|x\| \leq C_1\|A\|$. Therefore,
\begin{eqnarray}
\left(\sum_{j=1}^m \sum_{k=1}^m \left\|A(x_j^{(1)},
x_k^{(2)})\right\|^q\right)^{\frac{1}{q}} \!\!& = & \!\!
\left(\sum_{j=1}^m \left\|A_1^{(m)}(x_j^{(1)})
\right\|^q\right)^{\frac{1}{q}} \nonumber\\
\!\!& \leq & \!\! \pi_{q;r}(A_1^{(m)})\left
\|(x_j^{(1)})_{j=1}^m \right\|_{w,r} \nonumber\\
\!\!& \leq & \!\! C_2\|A_1^{(m)}\|\left \|(x_j^{(1)})_{j=1}^m
\right\|_{w,r}\nonumber\\
\!\!& = & \!\! C_2\sup_{x\in B_{E_1}}\left(\sum_{k=1}^m
\left\|A(x, x_k^{(2)})\right\|^q\right)^{\frac{1}{q}}\left
\|(x_j^{(1)})_{j=1}^m \right\|_{w,r}\nonumber\\
\!\!& = & \!\! C_2\sup_{x\in B_{E_1}}\left(\sum_{k=1}^m
\left\|A_x( x_k^{(2)})\right\|^q\right)^{\frac{1}{q}}\left
\|(x_j^{(1)})_{j=1}^m \right\|_{w,r}\nonumber\\
\!\!& \leq & \!\! C_2\sup_{x\in B_{E_1}}\pi_{q;p}(A_x)\left
\|(x_k^{(2)})_{j=1}^m \right\|_{w,p}\left \|(x_j^{(1)})_{j=1}^m
\right\|_{w,r}\nonumber\\
\!\!& \leq & \!\! C_1C_2\|A\|\left \|(x_j^{(1)})_{j=1}^m
\right\|_{w,r}\left \|(x_k^{(2)})_{j=1}^m \right\|_{w,p},\nonumber
\end{eqnarray}
which shows that $A$ is multiple $(q;r,p)$-summing and
$\pi_{q;r,p}(A) \leq  C_1C_2\|A\|$.\\
\indent Suppose now that the result holds for $n$, that is: for
every $E_1, \ldots, E_n \in B(p,q,r,F)$, ${\cal L}(E_1, \ldots,
E_n;F) = \Pi_{q;r, \ldots,r,p}^n(E_1, \ldots, E_n;F)$. To prove
the case $n+1$, let $E_1, \ldots, E_{n+1} \in B(p, q, r,F)$. As
$E_2, \ldots, E_{n+1}$ are $n$ Banach spaces in $B(p,q,r,F)$,
${\cal L}(E_2, \ldots, E_{n+1};F) = \Pi_{q;r, \ldots,r,p}^n(E_2,
\ldots, E_{n+1};F)$ by the induction hypotheses, so we can select
a constant $C_1$ such that
$$\pi_{q;r, \ldots,r,p}(B) \leq C_1\|B\| {\rm ~for~every~}B \in {\cal
L}(E_2, \ldots, E_{n+1};F).$$ Since $E_1 \in B(p,q,r,F)$, there is
a constant $C_2$ such that
$$\pi_{q;r}(v) \leq C_2\|v\| {\rm ~for~every~}v \in {\cal L}(E_1;\ell_q(F)).$$
Let $A \in {\cal L}(E_1, \ldots, E_{n+1};F)$. Given sequences
$(x_j^{(1)})_{j=1}^\infty \in \ell_r^w(E_1), \ldots,$
$(x_j^{(n)})_{j=1}^\infty \in \ell_r^w(E_n)$ and
$(x_j^{(n+1)})_{j=1}^\infty \in \ell_p^w(E_{n+1})$, fix $m \in
\mathbb{N}$ and consider the continuous linear operator
$$A_1^{(m)} \colon E_1 \longrightarrow \ell_q(F) ~:~ A_1^{(m)}(x) =
\left((A(x,x_{j_2}^{(2)}, \ldots,x_{j_{n+1}}^{(n+1)})_{j_2,
\ldots,j_{n+1}=1}^m,0,0, \ldots \right).
$$ So, $A_1^{(m)}$ is $(q;r)$-summing and $\pi_{q;r}(A_1^{(m)})
\leq C_2\|A_1^{(m)}\|$. For each $x \in B_{E_1}$, consider the
continuous $n$-linear mapping
$$A_x^n \colon E_2 \times \cdots \times E_{n+1} \longrightarrow F ~:~
A_x^n(x_2, \ldots,x_{n+1}) = A(x,x_2, \ldots, x_{n+1}). $$ So,
$\pi_{q;r, \ldots, r, p}(A_x^n) \leq C_1\|A_x^n\| \leq
C_1\|A\|\|x\| \leq C_1\|A\|$. Hence,
$$ \left(\sum_{j_1=1}^m \cdots \sum_{j_{n+1}=1}^m
\left\|A(x_{j_1}^{(1)}, \ldots
x_{j_{n+1}}^{(n+1)})\right\|^q\right)^{\frac{1}{q}}
=\left(\sum_{j=1}^m \left\|A_1^{(m)}(x_j^{(1)})
\right\|^q\right)^{\frac{1}{q}}\hspace*{60em} $$
$$ \leq
\pi_{q;r}(A_1^{(m)})\left \|(x_j^{(1)})_{j=1}^m \right\|_{w,r}
\leq C_2\|A_1^{(m)}\|\left \|(x_j^{(1)})_{j=1}^m
\right\|_{w,r}\hspace*{60em}$$
$$ = C_2\sup_{x\in B_{E_1}}\left(\sum_{j_2=1}^m \cdots \sum_{j_{n+1}=1}^m
\left\|A(x,x_{j_2}^{(2)},
\ldots,x_{j_{n+1}}^{(n+1)})\right\|^q\right)^{\frac{1}{q}}\left
\|(x_j^{(1)})_{j=1}^m \right\|_{w,r}\hspace*{60em}$$
$$ = C_2\sup_{x\in B_{E_1}}\left(\sum_{j_2=1}^m \cdots \sum_{j_{n+1}=1}^m \left\|
A_x^n( x_{j_2}^{(2)},
\ldots,x_{j_{n+1}}^{(n+1)}))\right\|^q\right)^{\frac{1}{q}}\left
\|(x_j^{(1)})_{j=1}^m \right\|_{w,r}\hspace*{60em}$$
$$ \leq C_2\sup_{x\in B_{E_1}}\pi_{q;r, \ldots, r,p}(A_x^n)
\left(\prod_{k=2}^n \left \|(x_j^{(k)})_{j=1}^m
\right\|_{w,r}\right)\left \|(x_j^{(n+1)})_{j=1}^m
\right\|_{w,p}\left \|(x_j^{(1)})_{j=1}^m
\right\|_{w,r}\hspace*{60em}$$
$$ \leq C_1C_2\|A\|\left(\prod_{k=1}^n \left \|(x_j^{(k)})_{j=1}^m
\right\|_{w,r}\right)\left \|(x_j^{(n+1)})_{j=1}^m
\right\|_{w,p},\hspace*{60em}$$ which shows that $A$ is multiple
$(q;r,\ldots, r,p)$-summing and completes the proof.
\end{proof}

Rewriting the proof above for $E_1 = \cdots = E_n = E$, we obtain

\begin{theorem} Let $p, r \in [1,q]$ and let
$E$ and $F$ be Banach spaces such that ${\cal L}(E;F) =
\Pi_{q;p}(E;F)$ and ${\cal L}(E;\ell_q(F)) =
\Pi_{q;r}(E;\ell_q(F))$ with
$$\pi_{q;p}(u) \leq C_1\|u\| {\it ~for~every~}u \in {\cal L}(E;F){\it ~and~}  $$
$$\pi_{q;r}(v) \leq C_2\|v\| {\it ~for~every~}v \in {\cal L}(E;\ell_q(F)).$$
Then, for every $n \geq 2$,
$${\cal L}(^nE;F) = \Pi_{q;r,\ldots, r,p}^n(^nE;F){\it~and} $$
$$\pi_{q;r, \ldots, r,p}(A) \leq C_1C_2^{n-1}\|A\| {\it ~for~every~}A \in {\cal L}(^nE;F).$$
\end{theorem}

For scalar-valued mappings we get the following particular cases:

\begin{corollary}\label{corolario} Given $1 \leq r \leq q$, by $B(r,q)$ we mean
the collection of all Banach spaces
$E$ such that ${\cal L}(E;\ell_q) = \Pi_{q;r}(E;\ell_q).$ Then,
for every $n \geq 2$,
$${\cal L}(E_1, \ldots, E_n) = \Pi_{q;r, \ldots,r,q}^n(E_1, \ldots,
E_n)$$ whenever $E_1, \ldots, E_n \in B(r,q)$.
\end{corollary}

\begin{corollary} Let $1 \leq r \leq q$ and let
$E$ be a Banach space such that ${\cal L}(E;\ell_q) =
\pi_{q;r}(E;\ell_q)$ with $\pi_{q;r}(v) \leq C\|v\|$ for every $v
\in {\cal L}(E;\ell_q)$. Then, for every $n \geq 2$,
$${\cal L}(^nE) = \Pi_{q;r,\ldots,r,q}^n(^nE){\it~and~} \pi_{q;r, \ldots, r,q}(A)
 \leq C^{n-1}\|A\| {\it ~for~every~}A \in {\cal L}(^nE).$$
\end{corollary}

\section{Applications}
We start by showing that some known coincidence theorems are easy
combinations of our general results with linear ones. From now on,
$n$ will always be an integer not smaller than 2.

\begin{proposition}{\rm \cite[Theorem 3.2]{bpgv}} If $F$ has cotype
$q$ and $E_1, \ldots, E_n$ are arbitrary
Banach spaces, then
$${\cal L}(E_1, \ldots, E_n;F) = \Pi_{q;1}^n(E_1, \ldots,
E_n;F){\it ~and~}$$
$$ \pi_{q;1}(A)
 \leq C_q(F)^n\|A\| {\it ~for~every~}A \in {\cal L}(E_1, \ldots, E_n;F),$$
where $C_q(F)$ is the cotype $q$ constant of $F$.
\end{proposition}

\begin{proof} Both $F$ and $\ell_q(F)$ have cotype $q$ (see \cite[Theorem
11.12]{djt}), so ${\cal L}(E;F) = \Pi_{q;1}(E;F)$ and ${\cal
L}(E;\ell_q(F)) = \Pi_{q;1}(E;\ell_q(F))$ for every Banach space
$E$ by \cite[Corollary 11.17]{djt}. The desired coincidence
follows from Theorem \ref{teorema} and the estimate for the norms
from its proof.
\end{proof}

\begin{proposition}\label{proposicao}{\rm \cite[Theorem 5.1]{bpgv}} If $E_1, \ldots,
E_n$ are ${\cal L}_1$-spaces and $H$ is a Hilbert space, then
$${\cal L}(E_1, \ldots, E_n;H) = \Pi_{2;2}^n(E_1, \ldots,
E_n;H){\it ~and~}$$
$$ \pi_{2;2}(A)
 \leq K_{G}^n\|A\| {\it ~for~every~}A \in {\cal L}(E_1, \ldots, E_n;H),$$
 where $K_{G}$ stands for the Grothendieck constant.
\end{proposition}

\begin{proof} Both $H$ and $\ell_2(H)$ are ${\cal L}_2$-spaces (see
\cite[Ex. 23.17(a)]{df}), so ${\cal L}(E;H) = \Pi_{2;2}(E;H)$ and
${\cal L}(E;\ell_2(H)) = \Pi_{2;2}(E;\ell_2(H))$ for every ${\cal
L}_1$-space $E$ by \cite[Theorems 3.1 and 2.8]{djt}. As before,
the result follows from Theorem \ref{teorema} and its proof.
\end{proof}

\begin{proposition} {\rm \cite[Theorem 3.1]{bpgv}} If $F$ has cotype
$2$ and $E_1, \ldots, E_n$ are ${\cal L}_\infty$-spaces, then
${\cal L}(E_1, \ldots, E_n;F) = \Pi_{2;2}^n(E_1, \ldots, E_n;F).$
\end{proposition}

\begin{proof} Both $F$ and $\ell_2(F)$ have cotype 2 \cite[Theorem 11.12]{djt},
 so ${\cal L}(E;F) = \Pi_{2;2}(E;F)$
and ${\cal L}(E;\ell_2(F)) = \Pi_{2;2}(E;\ell_2(F))$ for every
${\cal L}_\infty$-space $E$ by \cite[Theorem 11.14(a)]{djt}. Call
on Theorem \ref{teorema} once more.
\end{proof}

Using \cite[Theorem 11.14(b)]{djt} instead of \cite[Theorem
11.14(a)]{djt} in the proof above we obtain

\begin{proposition} If $F$ has cotype $q > 2$,
$E_1, \ldots, E_n$ are ${\cal L}_\infty$-spaces and $r < q$, then
${\cal L}(E_1, \ldots, E_n;F) = \Pi_{q;r}^n(E_1, \ldots, E_n;F).$
\end{proposition}

Now we derive some coincidence situations which, as far as we
know, are new. The first one complements nice information given in
\cite[Corollary 3.20]{arkiv}.

\begin{proposition} If $E_1, \ldots, E_n$ are arbitrary Banach spaces and $q \geq
2$, then
$${\cal L}(E_1, \ldots, E_n) = \Pi_{q;1,\ldots,1,q}^n(E_1, \ldots,
E_n){\it ~and~}$$
$$ \pi_{q;1,\ldots,1,q}(A)
 \leq C_q(\ell_q)^{n-1}\|A\| {\it ~for~every~}A \in {\cal L}(E_1, \ldots, E_n).$$
\end{proposition}

\begin{proof} Since $\ell_q$ has cotype $q$,
 ${\cal L}(E;\ell_q) = \Pi_{q;1}(E;\ell_q)$ for every Banach space
 $E$, so the desired coincidence follows from Corollary \ref{corolario} and the
 estimate for the norms follows from the proof of Theorem \ref{teorema}.
\end{proof}

In order to prove that there is no general inclusion theorem for
multiple summing multilinear mappings (that is, $q \leq p
\not\Longrightarrow \Pi_{q;q}^n \subseteq \Pi_{p;p}^n$), in
\cite[Theorem 3.6]{jmaa} the authors show that
$\Pi_{q;q}^2(^2\ell_1) \neq {\cal L}(^2\ell_1)$ for every $q>2$,
whereas $\Pi_{2;2}^2(^2\ell_1) = {\cal L}(^2\ell_1)$ (cf.
Proposition \ref{proposicao}). Next proposition shows that the
non-coincidence $\Pi_{q;q}^2(^2\ell_1) \neq {\cal L}(^2\ell_1)$,
$q > 2$, is quite sharp.

\begin{proposition} If $E_1, \ldots, E_n$ are ${\cal L}_1$-spaces and $2 \leq r < q$, then
${\cal L}(E_1, \ldots, E_n) = \Pi_{q;r,\ldots,r,q}^n(E_1, \ldots,
E_n).$
\end{proposition}

\begin{proof} For every ${\cal L}_1$-space $E$ and $2 \leq r < q$, ${\cal L}(E;\ell_q) =
 \Pi_{q;r}(E;\ell_q)$ by a result due to Bennet \cite[Proposition
5.2(iv)]{bennet}. The result follows from Corollary
\ref{corolario}.
\end{proof}

Using \cite[Proposition 5.1(ii)]{bennet} instead of
\cite[Proposition 5.2(iv)]{bennet} we get:

\begin{proposition} If $E_1, \ldots, E_n$ are ${\cal L}_\infty$-spaces
and $q > r, q >2$, then
${\cal L}(E_1, \ldots, E_n) = \Pi_{q;r,\ldots,r,q}^n(E_1, \ldots,
E_n).$
\end{proposition}

\section{Multiple summing mappings on ${\cal L}_1$-spaces}
We have already obtained some applications of our results to
multiple summing mappings on ${\cal L}_1$-spaces. In this section
we go a little further in this direction. Before using a new
approach, we apply our general results a couple of times more.

\begin{proposition} Let $E_1, \ldots, E_n$ be ${\cal L}_1$-spaces.\\
{\rm(a)} ${\cal L}(E_1, \ldots, E_n;H) =
\Pi_{q;1,\ldots,1,q}^n(E_1, \ldots, E_n;H)$ for every Hilbert space $H$ and any $q \geq 2$.\\
{\rm(b)} ${\cal L}(E_1, \ldots, E_n;F) =
\Pi_{q;1,\ldots,1,r}^n(E_1, \ldots, E_n;F)$ for every ${\cal
L}_q$-space $F$ and any $2 \leq r < q$.
\end{proposition}

\begin{proof} (a) Let $E$ be an ${\cal L}_1$-space and $H$ be a Hilbert space
. \cite[Theorems 3.1 and 2.8]{djt} yield that ${\cal L}(E; H) =
\Pi_{q;q}(E;H)$, and ${\cal L}(E;\ell_q(H)) =
\Pi_{q;1}(E;\ell_q(H))$ because $\ell_q(H)$ has
cotype $q$. Theorem \ref{teorema} gives the result.\\
\noindent (b) Let $E$ be an ${\cal L}_1$-space and $F$ be an
${\cal L}_q$-space. Using \cite[Proposition 5.2(iv)]{bennet} once
more we know that ${\cal L}(E;F) = \Pi_{q;r}(E;F)$. ${\cal L}(E;F)
= \Pi_{q;1}(E;\ell_q(F))$ as $\ell_q(F)$ has cotype $q$, so
 Theorem \ref{teorema} completes the proof.
\end{proof}

Next result allows us to go a little bit further.

\begin{theorem}\label{novoteorema} Let $r \geq s$. If ${\cal L}(\ell_1;F) =
\Pi_{r;s}(\ell_1;F)$, then $${\cal L}(E_1, \ldots, E_n;F) =
\Pi_{r;\,\min\{s,2\}}^n(E_1, \ldots, E_n;F)$$ for every $n \in
\mathbb{N}$ and any ${\cal L}_1$-spaces $E_1, \ldots, E_n$.
\end{theorem}

\begin{proof} By a standard localization argument we may assume
$E_1 = \cdots = E_n = \ell_1$. Let $(x_j^{(1)})_{j=1}^{m_1},
\ldots, (x_j^{(n)})_{j=1}^{m_n}$ be $n$ finite sequences in
$\ell_1$. \\
\noindent Claim: For every $1 \leq p \leq 2$,
$$\left \| (x_{j_1}^{(1)} \otimes \cdots \otimes x_{j_n}^{(n)})
_{j_1, \ldots,  j_n = 1}^{m_1, \ldots, m_n}\right \|_{w,p} \leq
K_G^{2n-2}
 \left \| (x_j^{(1)})_{j=1}^{m_1}\right \|_{w,p} \cdots \left \|
 (x_j^{(n)})_{j=1}^{m_n}\right \|_{w,p}$$
on the completed $n$-fold projective tensor product
$\widehat\otimes_\pi^n\ell_1$.\\
\noindent Proof of the claim: we proceed by induction on $n$.
Given $A \in {\cal L}(^2\ell_1)$, by \cite[Theorem 3.4]{jmaa} we
know that $A \in \Pi_{p;p}^2(^2\ell_1)$ and $\pi_{p;p}(A) \leq
K_G^2 \|A\|$. Denoting by $\sigma$ the canonical bilinear mapping
from $\ell_1 \times \ell_1$ to $\ell_1 \hat\otimes_\pi \ell_1$,
$\sigma(x,y) = x \otimes y$, and taking the supremum over all
$\varphi \in B_{(\ell_1 \widehat\otimes_\pi \ell_1)'}$, it follows
that
\begin{eqnarray}
\sup_{\varphi} \left ( \sum_{j_1, j_2=1}^{m_1,m_2}\left
|\varphi(x_j^{(1)} \otimes x_j^{(2)}) \right |^p\right
)^\frac{1}{p} \!\! & = & \!\! \sup_{\varphi} \left ( \sum_{j_1,
j_2=1}^{m_1,m_2}\left
|\varphi \circ \sigma(x_j^{(1)},x_j^{(2)}) \right |^p\right )^\frac{1}{p}\nonumber\\
\!\! & \leq & \!\! \sup_{\varphi} \pi_{p;p}(\varphi \circ \sigma)
\left \| (x_j^{(1)})_{j=1}^{m_1}\right \|_{w,p}\left \|
(x_j^{(2)})_{j=1}^{m_2}\right \|_{w,p}\nonumber\\
\!\! & \leq & \!\! K_G^2\sup_{\varphi} \|\varphi \circ \sigma\|
\left \| (x_j^{(1)})_{j=1}^{m_1}\right \|_{w,p}\left \|
(x_j^{(2)})_{j=1}^{m_2}\right \|_{w,p}\nonumber\\
\!\! & = & \!\! K_G^2 \left \| (x_j^{(1)})_{j=1}^{m_1}\right
\|_{w,p}\left \| (x_j^{(2)})_{j=1}^{m_2}\right \|_{w,p}.
\end{eqnarray}
Now suppose that the desired inequality holds for $k$ and let us
prove that it holds for $k+1$. We are assuming that
$$\left \| (x_{j_1}^{(1)} \otimes \cdots \otimes x_{j_k}^{(k)})_{j_1, \ldots,
 j_k= 1}^{m_1, \ldots, m_k}\right \|_{w,p} \leq K_G^{2k-2}
 \left \| (x_j^{(1)})_{j=1}^{m_1}\right \|_{w,p} \cdots \left \|
 (x_j^{(k)})_{j=1}^{m_k}\right \|_{w,p}.$$
Using that $\widehat\otimes_\pi^k\ell_1$ is isometrically
isomorphic to $\ell_1$, applying first (1) then the induction
hypotheses, we get
$$\left \| (x_{j_1}^{(1)} \otimes \cdots \otimes
x_{j_k}^{(k)}\otimes x_{j_{k+1}}^{(k+1)})_{j_1, \ldots,
 j_{k+1}= 1}^{m_1, \ldots, m_{k+1}}\right \|_{w,p} \hspace*{25em}$$
$$\hspace*{3em} = \left \| \left(\left(x_{j_1}^{(1)} \otimes \cdots \otimes
x_{j_k}^{(k)}\right)\otimes x_{j_{k+1}}^{(k+1)} \right)_{j_1,
\ldots,
 j_{k+1}= 1}^{m_1, \ldots, m_{k+1}}\right \|_{w,p} \hspace*{25em}$$
 $$\hspace*{3em} \leq K_G^2
 \left \| \left(x_{j_1}^{(1)} \otimes \cdots \otimes
x_{j_k}^{(k)}\right)_{j_1, \ldots,
 j_{k}= 1}^{m_1, \ldots, m_{k}}\right \|_{w,p} \cdot \left \|
 (x_{j}^{(k+1)})_{j=1}^{m_{k+1}}\right \|_{w,p}\hspace*{25em}$$
$$\hspace*{3em} \leq K_G^2 K_G^{2k-2}
 \left \| (x_j^{(1)})_{j=1}^{m_1}\right \|_{w,p} \cdots \left \|
 (x_j^{(k)})_{j=1}^{m_k}\right \|_{w,p}\left \|
 (x_{j}^{(k+1)})_{j=1}^{m_{k+1}}\right \|_{w,p}\hspace*{25em}$$
$$\hspace*{3em} = K_G^{2(k+1) -2}
 \left \| (x_j^{(1)})_{j=1}^{m_1}\right \|_{w,p} \cdots \left \|
 (x_j^{(k)})_{j=1}^{m_k}\right \|_{w,p}\left \|
 (x_{j}^{(k+1)})_{j=1}^{m_{k+1}}\right \|_{w,p},\hspace*{25em}$$
completing the proof of the claim.\\
\indent Let $A \in {\cal L}(^n\ell_1;F)$. By $A_L$ we mean the
linearization of $A$ on $\widehat\otimes_\pi^n\ell_1$, that is
$A_L \in {\cal L}(\widehat\otimes_\pi^n\ell_1;F)$ and $A_L(x_1
\otimes \cdots \otimes x_n) = A(x_1, \ldots, x_n)$ for every $x_j
\in \ell_1$. Since $\widehat\otimes_\pi^n\ell_1$ is isometrically
isomorphic to $\ell_1$, by assumption we have that $A_L$ is
$(r;s)$-summing and $\pi_{r;s}(A_L) \leq M\|A_L\| = M\|A\|$, where
$M$ is a constant independent of $A$. Using the claim with $p
=\min\{s,2\}$ we get
$$\left(\sum_{j_1, \ldots, j_n= 1}^{m_1,\ldots, m_n} \left\|A(x_{j_1}^{(1)},
\ldots,
x_{j_n}^{(n)})\right\|^r\right)^\frac{1}{r}\hspace*{25em}$$
$$ \hspace*{3em} = \left(\sum_{j_1, \ldots, j_n= 1}^{m_1,\ldots, m_n} \left\|A_L(x_{j_1}^{(1)}
\otimes \cdots \otimes x_{j_n}^{(n)})\right\|^r\right)^\frac{1}{r}
\hspace*{25em}$$
$$\hspace*{3em} \leq \pi_{r;s}(A_L)  \left \| (x_{j_1}^{(1)} \otimes \cdots \otimes x_{j_n}^{(n)})
_{j_1, \ldots,  j_n= 1}^{m_1, \ldots, m_n}\right \|_{w,s}
\hspace*{25em}$$
$$\hspace*{3em} \leq \pi_{r;s}(A_L)  \left \| (x_{j_1}^{(1)} \otimes \cdots \otimes x_{j_n}^{(n)})
_{j_1, \ldots,  j_n= 1}^{m_1, \ldots, m_n}\right
\|_{w,\min\{s,2\}} \hspace*{25em}$$
$$\hspace*{3em} \leq M \|A\|
K_G^{2n-2}
 \left \| (x_j^{(1)})_{j=1}^{m_1}\right \|_{w,\min\{s,2\}} \cdots \left \|
 (x_j^{(n)})_{j=1}^{m_n}\right \|_{w,\min\{s,2\}}, \hspace*{25em}$$
which shows that $A$ is multiple $(r;\min\{s,2\})$-summing.
\end{proof}

\begin{corollary} Let $1 \leq p \leq 2$, $r \geq p$ and let $F$ be a
Banach space. The following assertions are equivalent:\\
{\rm (a)} ${\cal L}(\ell_1;F) =
\Pi_{r;p}(\ell_1;F)$.\\
{\rm (b)} ${\cal L}(^n\ell_1;F) =
\Pi_{r;\,p}^n(^n\ell_1;F)$ for every $n \in \mathbb{N}$.\\
{\rm (c)} ${\cal L}(^n\ell_1;F) = \Pi_{r;\,p}^n(^n\ell_1;F)$ for
some $n \in \mathbb{N}$.
\end{corollary}

\begin{proof} (a) $\Longrightarrow$ (b) follows from Theorem
\ref{novoteorema}, (b) $\Longrightarrow$ (c) is obvious and (c)
$\Longrightarrow$ (a) follows from Remark \ref{remark}.
\end{proof}

The last result of the paper makes clear how our methods
systematize the subject and generalize and simplify the known
results. The coincidence
\begin{eqnarray} {\cal L}(^n\ell_1;\ell_2) = \Pi_{p;p}^n(^n\ell_1;\ell_2)
{\rm ~for~}n \in \mathbb{N} {\rm ~and~}1 \leq p \leq 2
\end{eqnarray}
was proved in \cite{davidtese} in the following fashion: the
author proves first that ${\cal L}(^n\ell_1;\ell_2) =
\Pi_{2;2}^n(^n\ell_1;\ell_2)$ (a result we reobtained in
Proposition \ref{proposicao}), then uses this to prove that ${\cal
L}(^n\ell_1;\ell_2) = \Pi_{1;1}^n(^n\ell_1;\ell_2)$ (see also
\cite[Theorem 5.2]{bpgv}), and finally uses this last coincidence
to obtain (2) (cf. \cite[Corolario 5.24]{davidtese}). Several
auxiliary results are used along the way. On the other hand, (2)
is nothing but a particular case of the next corollary, which is a
straightforward combination of Theorem \ref{novoteorema} with
\cite[Theorem 5.2]{bennet}.

\begin{corollary} Let $E_1, \ldots, E_n$ be ${\cal L}_1$-spaces
and let $F$ be an ${\cal L}_q$-space, $1 \leq q < +\infty$. Then
$${\cal L}(E_1 \ldots, E_n;F) = \Pi_{r;p}^n(E_1 \ldots, E_n;F)$$
if either\\
\indent {\rm (a)} $q < 2$, $r \geq q^*$ and $p =2$, where
$\frac{1}{q}
 + \frac{1}{q^*} = 1$, or\\
\indent {\rm (b)} $q >2$, $r \geq q$ and $p = 2$, or\\
\indent {\rm (c)} $q = 2$ and $1 \leq p \leq r \leq 2$.
\end{corollary}

\medskip
\noindent {\sc Acknowledgement.} Part of this paper was written
while G.B. was a CNPq Postdoctoral Fellow in the Departamento de
An\'alisis Matem\'atico at Universidad de Valencia. He thanks
Pilar Rueda and the members of the department for their kind
hospitality.

\vspace*{1em}
 \noindent [Geraldo Botelho] Faculdade de Matem\'atica, Universidade Federal de Uberl\^andia,
38.400-902 - Uberl\^andia, Brazil,  e-mail: botelho@ufu.br.

\medskip

\noindent [Daniel Pellegrino] Departamento de Matem\'atica,
Universidade Federal da Pa-ra\'iba, 58.051-900 - Jo\~ao Pessoa,
Brazil, e-mail: dmpellegrino@gmail.com.

\end{document}